\documentclass[11pt]{amsart}
\usepackage{amsmath,amssymb,amsthm}
\usepackage[margin=1.15in]{geometry}
\usepackage[hidelinks]{hyperref}

\makeatletter
\@namedef{subjclassname@2020}{\textup{2020} Mathematics Subject Classification}
\makeatother

\newtheorem{theorem}{Theorem}[section]
\newtheorem{lemma}[theorem]{Lemma}
\newtheorem{proposition}[theorem]{Proposition}
\newtheorem{corollary}[theorem]{Corollary}

\newtheorem*{theoremA}{Theorem A}
\newtheorem*{theoremB}{Theorem B}
\newtheorem*{theoremC}{Theorem C}
\theoremstyle{definition}
\newtheorem{definition}[theorem]{Definition}
\newtheorem{example}[theorem]{Example}

\theoremstyle{remark}
\newtheorem{remark}[theorem]{Remark}

\newcommand{\R}{\mathbb{R}}
\newcommand{\N}{\mathbb{N}}
\newcommand{\symd}{\mathbin{\triangle}}
\newcommand{\eps}{\varepsilon}
\newcommand{\Leb}{\lambda}

\begin{document}

\title{Open segments on topological measure spaces}

\author{Mohammad Javaheri}
\address{Department of Mathematics and Statistics, Siena University\\
515 Loudon Road, Loudonville, NY 12211, USA}
\email{mjavaheri@siena.edu}

\date{November 2, 2026}

\subjclass[2020]{Primary 28A25; Secondary 28C15, 46G10}
\keywords{Non-atomic measure, segment, open segment, Lyapunov convexity theorem, weakly outer regular measure}

\begin{abstract}
Let $\mu$ be a finite non-atomic measure on a set $X$. A \emph{segment} for
$\mu$ is an increasing family $u(t)$, $t\in[0,1]$, of measurable sets with
$u(0)=\emptyset$, $u(1)=X$ and $\mu(u(t))=t\,\mu(X)$; by a theorem of
Fryszkowski, finitely many measures $f\,d\mu$ admit a common segment. We study
\emph{open segments}, i.e.\ segments consisting of open sets, when $X$ carries a
topology. We show that weak outer regularity of $\mu$ is equivalent to the
approximability of measurable sets by open sets, and that it holds for every
finite Borel measure on a metrizable space. Our main result states that if
$\mu$ is finite, non-atomic and weakly outer regular, then for all
$f,g\in L^1(\mu)$ with $(\int f\,d\mu,\int g\,d\mu)\neq(0,0)$ there is a common
open segment for $f\,d\mu$ and $g\,d\mu$; in particular, for every
$f\in L^1(\mu)$ there is an open segment for $\mu$ which is simultaneously a
segment for $f\,d\mu$. Relative versions allow one to prescribe an initial set.
We also characterize open segments as upper semicontinuous functions with
uniform distribution and show that on separable metrizable spaces a single measure admits such a function which is continuous.
\end{abstract}

\maketitle

%======================================================================
\section{Introduction}
%======================================================================

A measure space $(X,\Sigma,\mu)$ is \emph{non-atomic} if every $A\in\Sigma$
with $\mu(A)>0$ contains some $B\in\Sigma$ with $0<\mu(B)<\mu(A)$.
Sierpi\'nski~\cite{Si} proved that the range of a finite non-atomic measure
$\mu$ is the whole interval $[0,\mu(X)]$, and Lyapunov~\cite{Ly} proved that the
range of a finite-dimensional non-atomic vector measure is compact and convex;
see \cite{Ha,Li,HR} for other proofs. There is an \emph{ordered} refinement of
these facts, in which the intermediate values are realized along a single
increasing family of sets.

\begin{definition}\label{def:segment}
Let $(X,\Sigma,\mu)$ be a finite measure space and let $\eta$ be a finite
signed measure on $\Sigma$. A \emph{segment} for $\eta$ is a map
$u\colon[0,1]\to\Sigma$ such that $u(0)=\emptyset$, $u(1)=X$, $u(s)\subseteq
u(t)$ whenever $s\le t$, and
\begin{equation}\label{eq:segment}
\eta(u(t))=t\,\eta(X)\qquad\text{for all } t\in[0,1].
\end{equation}
For $f\in L^1(X,\Sigma,\mu)$ we write $\mu_f(E)=\int_E f\,d\mu$, and a segment
for $\mu_f$ is also called a segment for $f$. If $X$ is a topological space and
$u(t)$ is open for every $t$, then $u$ is an \emph{open segment}, or
\emph{O-segment}.
\end{definition}

\begin{theorem}[Fryszkowski \cite{Fr1,Fr2}; cf.\ Cellina--Colombo--Fonda \cite{CCF}]\label{thm:Fr-intro}
Let $(X,\Sigma,\mu)$ be a finite non-atomic measure space and let
$S\subseteq L^1(\mu)$ be finite. Then there is a segment for $\mu$ which is a
segment for every $f\in S$.
\end{theorem}

We give an elementary proof in Section~\ref{sec:measurable}, together with the
facts that finiteness of $S$ cannot be dropped and that an analogue holds on
$\sigma$-finite spaces.

In this paper $X$ is in addition a topological space whose open sets are
measurable, and we ask when the sets $u(t)$ can be chosen open. Some
hypothesis linking $\mu$ to the topology is needed: in
Example~\ref{ex:no-O-segment} the values of $\mu$ on open sets fill
$[0,\mu(X)]$, yet $\mu$ has no O-segment. The natural hypothesis is the
following.

\begin{definition}\label{def:regularity}
Let $\mu$ be a finite measure on a $\sigma$-algebra $\Sigma$ of subsets of a
topological space $X$ that contains all open sets, and let $A\in\Sigma$.
\begin{enumerate}
\item[(i)] $A$ is \emph{approximable by open sets} if for every $\eps>0$ there
is an open set $O$ with $\mu(A\symd O)<\eps$.
\item[(ii)] $\mu$ is \emph{weakly outer regular at $A$} if for every $\eps>0$
there is an open set $O$ with $\mu(A\setminus O)=0$ and $\mu(O\setminus
A)<\eps$.
\item[(iii)] $\mu$ is \emph{outer regular at $A$} if for every $\eps>0$ there is
an open set $O\supseteq A$ with $\mu(O\setminus A)<\eps$.
\end{enumerate}
We say that $\mu$ is \emph{weakly outer regular} if (ii) holds for all
$A\in\Sigma$.
\end{definition}

Clearly (iii)$\Rightarrow$(ii)$\Rightarrow$(i). Our first observation is that,
as properties of the measure, (i) and (ii) coincide, and that they are
automatic in the most common setting.

\begin{theoremA}[Lemma~\ref{lem:approx-equiv}, Proposition~\ref{prop:Gdelta}]
A finite measure $\mu$ is weakly outer regular if and only if every measurable
set is approximable by open sets. If every closed subset of $X$ is a $G_\delta$
set (e.g.\ if $X$ is metrizable), then every finite measure on the Borel
$\sigma$-algebra of $X$, or on its completion, is outer regular, hence weakly
outer regular.
\end{theoremA}

The main result of the paper concerns common O-segments.

\begin{theoremB}[Theorem~\ref{thm:main}]
Let $X$ be a topological space and let $\mu$ be a finite, non-atomic, weakly
outer regular measure on a $\sigma$-algebra containing the open sets.
\begin{enumerate}
\item[(a)] For every $f\in L^1(\mu)$ there is an O-segment for $\mu$ which is
also a segment for $f$.
\item[(b)] For all $f,g\in L^1(\mu)$ with
$(\mu_f(X),\mu_g(X))\neq(0,0)$ there is a common O-segment for $f$ and $g$.
\end{enumerate}
Moreover, the O-segment may be required to pass through a prescribed open set
$D$ (Theorem~\ref{thm:main}(c)), and there is a relative version starting
from an arbitrary measurable set (Theorem~\ref{thm:relative}).
\end{theoremB}

If $\mu_f(X)=\mu_g(X)=0$, the condition~\eqref{eq:segment} only says
$\mu_f(u(t))=\mu_g(u(t))=0$, and the family $u(t)=\emptyset$ for $t<1$,
$u(1)=X$, is a (degenerate) common O-segment. The meaningful version of (b) in
that case asks in addition for a segment for $\mu$, which is a problem about
three functions $1,f,g$.

Finally, O-segments have a functional description, which gives a convenient
language for the open problem and a short treatment of metrizable spaces.

\begin{theoremC}[Propositions~\ref{prop:usc} and \ref{prop:continuous}]
An O-segment for $\mu$ that is also a segment for each $f\in S$ is the same as
(the family of strict sublevel sets of) an upper semicontinuous function
$\varphi\colon X\to[0,1]$ with $\varphi_*\mu_f=\mu_f(X)\,\Leb$ for $f\in
S\cup\{1\}$, where $\Leb$ is Lebesgue measure on $[0,1]$. If $X$ is separable
and metrizable and $\mu$ is a finite non-atomic Borel measure, then there is a
\emph{continuous} $\varphi\colon X\to[0,1]$ with $\varphi_*\mu=\mu(X)\Leb$.
\end{theoremC}

For compact manifolds and measures of full support vanishing on the boundary,
the last statement also follows from the homeomorphic measures theorem of
Oxtoby and Ulam~\cite{OU} (see also~\cite{AP}). Theorem~C shows that for a
single measure on a metrizable space the existence of O-segments is easy;
the content of Theorem~B lies in handling two densities at once, where no such
shortcut is available. We conjecture that Theorem~B(a) holds for any finite
number of functions.

\subsection*{Organization}
Section~\ref{sec:prelim} contains preliminaries, Theorem~A and the functional
description of O-segments. Section~\ref{sec:measurable} treats measurable
segments. Section~\ref{sec:single} treats O-segments for one measure, and
Section~\ref{sec:pairs} proves Theorem~B. 

%======================================================================
\section{Preliminaries}\label{sec:prelim}
%======================================================================

Throughout, $(X,\Sigma,\mu)$ is a measure space and $L^1=L^1(X,\Sigma,\mu)$
consists of real-valued integrable functions. For $f\in L^1$ the set function
$\mu_f$ is a finite signed measure. We write $A\symd B$ for the symmetric
difference. Whenever a topology on $X$ is present, we assume that all open sets
belong to $\Sigma$ and call $(X,\Sigma,\mu)$ a \emph{topological measure
space}.

\subsection{Segments}

The following elementary observations are used repeatedly.

\begin{remark}\label{rem:linear}
(1) If $u$ is a segment for the signed measures $\eta_1,\eta_2$, it is a
segment for every linear combination $a\eta_1+b\eta_2$. Hence the set of
$h\in L^1$ for which $u$ is a segment is a linear subspace.

(2) If $\eta(X)=0$, condition~\eqref{eq:segment} reads $\eta(u(t))=0$ and is
satisfied by the degenerate family $u(t)=\emptyset$ ($t<1$), $u(1)=X$. This is
why, for families of functions, we ask for segments that are \emph{also}
segments for $\mu$.

(3) Let $\mu(X)>0$ and $f\in L^1$. Put $f^*=\mu(X)f/\mu_f(X)$ if
$\mu_f(X)\neq0$, and $f^*=f+1$ if $\mu_f(X)=0$. Then $\mu_{f^*}(X)=\mu(X)$,
and every segment for both $\mu$ and $f^*$ is a segment for $f$, by~(1).
\end{remark}

\subsection{Weak outer regularity}

\begin{lemma}\label{lem:approx-equiv}
Let $\mu$ be a finite measure on a topological measure space and let
$A\in\Sigma$. Then $A$ is approximable by open sets if and only if $\mu$ is
weakly outer regular at $A$. In particular, $\mu$ is weakly outer regular if
and only if every measurable set is approximable by open sets.
\end{lemma}

\begin{proof}
If $\mu(A\setminus O)=0$ and $\mu(O\setminus A)<\eps$, then $\mu(A\symd
O)<\eps$. Conversely, suppose $A$ is approximable and let $\eps>0$. Choose open
sets $O_k$ with $\mu(A\symd O_k)<\eps 2^{-k}$, $k\ge1$, and put
$O=\bigcup_k O_k$. Then $\mu(A\setminus O)\le\mu(A\setminus O_k)<\eps2^{-k}$
for all $k$, so $\mu(A\setminus O)=0$, and $\mu(O\setminus
A)\le\sum_k\mu(O_k\setminus A)<\eps$.
\end{proof}

\begin{lemma}\label{lem:relative-wor}
Let $\mu$ be finite and weakly outer regular, let $Y$ be open and let
$A\in\Sigma$ with $A\subseteq Y$. For every $\eps>0$ there is an open set
$W\subseteq Y$ with $\mu(A\setminus W)=0$ and $\mu(W\setminus A)<\eps$.
\end{lemma}

\begin{proof}
If $W_0$ is open with $\mu(A\setminus W_0)=0$ and $\mu(W_0\setminus A)<\eps$,
then $W=W_0\cap Y$ is open, $A\setminus W=A\setminus W_0$ and $W\setminus
A\subseteq W_0\setminus A$.
\end{proof}

\begin{proposition}\label{prop:Gdelta}
Let $X$ be a topological space in which every closed set is a $G_\delta$ set,
for instance a metrizable space, and let $\mu$ be a finite measure on the Borel
$\sigma$-algebra $\mathcal{B}(X)$ or on its $\mu$-completion. Then $\mu$ is
outer regular at every measurable set; in particular $\mu$ is weakly outer
regular.
\end{proposition}

\begin{proof}
This is the standard argument for metric spaces, cf.~\cite[Theorem~1.1]{Bi}.
Let $\mathcal{A}$ be the family of $A\in\mathcal{B}(X)$ such that for every
$\eps>0$ there are a closed $F$ and an open $O$ with $F\subseteq A\subseteq O$
and $\mu(O\setminus F)<\eps$. If $F$ is closed, write $F=\bigcap_n O_n$ with
$O_n$ open and decreasing; then $\mu(O_n\setminus F)\to0$, so
$F\in\mathcal{A}$. The family $\mathcal{A}$ is closed under complements (swap
$F$ and $O$). If $A_k\in\mathcal{A}$ and $F_k\subseteq A_k\subseteq O_k$ with
$\mu(O_k\setminus F_k)<\eps2^{-k-1}$, put $O=\bigcup_kO_k$ and choose $N$ with
$\mu(\bigcup_k F_k\setminus\bigcup_{k\le N}F_k)<\eps/2$; then
$F=\bigcup_{k\le N}F_k$ is closed, $F\subseteq\bigcup_kA_k\subseteq O$ and
$\mu(O\setminus F)<\eps$. Thus $\mathcal{A}$ is a $\sigma$-algebra containing
the closed sets, so $\mathcal{A}=\mathcal{B}(X)$. Finally, a set in the
completion has the form $B\cup N$ with $B$ Borel and $N\subseteq N'$ for a
Borel null set $N'$; if $O\supseteq B$ and $O'\supseteq N'$ are open with
$\mu(O\setminus B)<\eps/2$ and $\mu(O')<\eps/2$, then $O\cup O'\supseteq B\cup
N$ and $\mu((O\cup O')\setminus(B\cup N))<\eps$.
\end{proof}

We also use the absolute continuity of the integral.

\begin{lemma}[{see e.g.\ \cite{Fo}}]\label{lem:abscont}
For $h\in L^1$ and $\eps>0$ there is $\delta>0$ such that $\int_E|h|\,d\mu<\eps$
whenever $E\in\Sigma$ and $\mu(E)<\delta$.
\end{lemma}

\subsection{O-segments as semicontinuous functions}

For a finite signed measure $\eta$ on $\Sigma$ and a measurable
$\varphi\colon X\to[0,1]$, let $\varphi_*\eta(E)=\eta(\varphi^{-1}(E))$ for
Borel $E\subseteq[0,1]$, and let $\Leb$ denote Lebesgue measure on $[0,1]$.
Recall that $\varphi$ is upper semicontinuous if and only if
$\{\varphi<t\}$ is open for every $t\in\R$.

\begin{proposition}\label{prop:usc}
Let $(X,\Sigma,\mu)$ be a finite topological measure space and
$S\subseteq L^1$, and put $S_1=S\cup\{1\}$.
\begin{enumerate}
\item[(a)] Let $u$ be an O-segment for $\mu$ which is a segment for every
$f\in S$, and let $\varphi(x)=\inf\{t\in[0,1]:x\in u(t)\}$. Then $\varphi$ is
upper semicontinuous, $\{\varphi<t\}=\bigcup_{s<t}u(s)$ for $t\in(0,1]$, and
$\varphi_*\mu_f=\mu_f(X)\,\Leb$ for all $f\in S_1$.
\item[(b)] Conversely, if $\varphi\colon X\to[0,1]$ is upper semicontinuous
and $\varphi_*\mu_f=\mu_f(X)\,\Leb$ for all $f\in S_1$, then $u(t)=\{\varphi<t\}$
for $t<1$ and $u(1)=X$ define an O-segment for $\mu$ which is a segment for
every $f\in S$.
\end{enumerate}
\end{proposition}

\begin{proof}
(a) Since $u(1)=X$, $\varphi$ takes values in $[0,1]$. If $\varphi(x)<t$, then
$x\in u(s)$ for some $s<t$; conversely $x\in u(s)$ with $s<t$ gives
$\varphi(x)\le s<t$. Hence $\{\varphi<t\}=\bigcup_{s<t}u(s)=\bigcup_n
u(t-1/n)$ is open for $t\in(0,1]$, while $\{\varphi<t\}=\emptyset$ for $t\le0$
and $=X$ for $t>1$; thus $\varphi$ is upper semicontinuous and measurable. For
$f\in S_1$ and $t\in(0,1]$, countable additivity gives
$\mu_f(\{\varphi<t\})=\lim_n(t-1/n)\mu_f(X)=t\,\mu_f(X)$. The signed measures
$\varphi_*\mu_f$ and $\mu_f(X)\Leb$ thus agree on the intervals $[0,t)$ and on
$[0,1]$; these form a $\pi$-system generating the Borel sets of $[0,1]$, so
the two signed measures coincide by Dynkin's $\pi$--$\lambda$ theorem applied
to their difference.

(b) The sets $u(t)$ are open and increasing, $u(0)=\{\varphi<0\}=\emptyset$, and
$\mu_f(u(t))=\varphi_*\mu_f([0,t))=t\,\mu_f(X)$ for $t<1$ and $f\in S_1$.
\end{proof}

Thus the existence of an O-segment for $\mu$ that is a segment for every $f\in
S$ is equivalent to the existence of an upper semicontinuous function that
pushes each $\mu_f$, $f\in S_1$, to the appropriate multiple of Lebesgue
measure. The measurable analogue (Theorem~\ref{thm:Fr}) is the existence of
such a measurable function.

%======================================================================
\section{Measurable segments}\label{sec:measurable}
%======================================================================

\subsection{Finite families}

\begin{lemma}\label{lem:halving}
Let $F\colon[0,1]\to\R$ be continuous with $F(0)=0$ and $F(1)=1$. Then there
is $c\in[0,1/2]$ with $F(c+1/2)-F(c)=1/2$.
\end{lemma}

\begin{proof}
The continuous function $G(c)=F(c+1/2)-F(c)-1/2$ on $[0,1/2]$ satisfies
$G(0)+G(1/2)=F(1)-F(0)-1=0$, so it has a zero by the intermediate value
theorem.
\end{proof}

\begin{theorem}\label{thm:Fr}
Let $(X,\Sigma,\mu)$ be a finite non-atomic measure space and let
$S\subseteq L^1$ be finite. Then there is a segment for $\mu$ which is a
segment for every $f\in S$.
\end{theorem}

\begin{proof}
We may assume $\mu(X)>0$, and by Remark~\ref{rem:linear}(3) that
$\mu_f(X)=\mu(X)$ for all $f\in S$. We prove by induction on $n\ge0$ the
statement:

\smallskip
\noindent$(\mathrm{S}_n)$ \emph{If $(X,\Sigma,\mu)$ is finite and non-atomic
with $\mu(X)>0$, and $f_1,\dots,f_n\in L^1$ satisfy $\mu_{f_i}(X)=\mu(X)$,
then there is a segment for $\mu$ which is a segment for each $f_i$.}
\smallskip

Fix such data and call a pair $(u,\Lambda)$ a \emph{partial segment} if
$\{0,1\}\subseteq\Lambda\subseteq[0,1]$, $u\colon\Lambda\to\Sigma$,
$u(0)=\emptyset$, $u(1)=X$, $u(s)\subseteq u(t)$ for $s\le t$ in $\Lambda$, and
$\mu(u(t))=\mu_{f_i}(u(t))=t\mu(X)$ for all $t\in\Lambda$ and $i\le n$. Order
partial segments by extension. The union of a chain is a partial segment, so
by Zorn's lemma there is a maximal partial segment $(u,\Lambda)$.

\emph{$\Lambda$ is closed.} Let $t'\in\overline\Lambda\setminus\Lambda$. There
is a monotone sequence $t_k\in\Lambda$ with $t_k\to t'$; let $u(t')$ be
$\bigcup_ku(t_k)$ if $(t_k)$ is increasing and $\bigcap_ku(t_k)$ if it is
decreasing. For $s\in\Lambda$, $s<t'$ implies $u(s)\subseteq u(t')$, and
$s>t'$ implies $u(t')\subseteq u(s)$; by countable additivity
$\mu(u(t'))=\mu_{f_i}(u(t'))=t'\mu(X)$. This extends $(u,\Lambda)$,
contradicting maximality.

\emph{$\Lambda=[0,1]$.} Otherwise $[0,1]\setminus\Lambda$ is a non-empty
relatively open set, so it contains a component $(a,b)$ with $a,b\in\Lambda$.
Let $X_0=u(b)\setminus u(a)$; then $\mu(X_0)=\mu_{f_i}(X_0)=(b-a)\mu(X)>0$. It
suffices to find $Z\in\Sigma$ with $u(a)\subseteq Z\subseteq u(b)$ and
$\mu(Z)=\mu_{f_i}(Z)=s\mu(X)$ for some $s\in(a,b)$ and all $i$, since then
$u(s)=Z$ extends $(u,\Lambda)$.

If $n=0$, choose $E\subseteq X_0$ with $0<\mu(E)<\mu(X_0)$ (possible since
$\mu$ is non-atomic) and let $Z=u(a)\cup E$. If $n\ge1$, apply
$(\mathrm{S}_{n-1})$ to $X_0$ with the restricted measure and
$f_1,\dots,f_{n-1}$; this gives a segment $v$ on $X_0$ with
$\mu(v(t))=\mu_{f_i}(v(t))=t\mu(X_0)$ for $i<n$. The function
$F(t)=\mu_{f_n}(v(t))/\mu(X_0)$ satisfies $F(0)=0$, $F(1)=1$, and it is
continuous: $|F(t)-F(s)|\le\mu(X_0)^{-1}\int_{v(t)\symd v(s)}|f_n|\,d\mu$ and
$\mu(v(t)\symd v(s))=|t-s|\mu(X_0)$, so Lemma~\ref{lem:abscont} applies. By
Lemma~\ref{lem:halving} there is $c$ with $F(c+1/2)-F(c)=1/2$. Put
$Z=u(a)\cup\bigl(v(c+1/2)\setminus v(c)\bigr)$. Then
$\mu(Z)=\mu_{f_i}(Z)=a\mu(X)+\tfrac12\mu(X_0)=\tfrac{a+b}2\mu(X)$ for all
$i\le n$.
\end{proof}

\begin{remark}\label{rem:nonatomic-iff}
Conversely, if a finite measure $\mu$ admits a segment $u$, then $\mu$ is
non-atomic: for $Y\in\Sigma$ with $\mu(Y)>0$, the function $t\mapsto\mu(Y\cap
u(t))$ is non-decreasing and Lipschitz with constant $\mu(X)$, with
values $0$ and $\mu(Y)$ at the endpoints, so it takes the value $\mu(Y)/2$, and
$Y$ is not an atom.
\end{remark}

\subsection{Infinite families}

\begin{proposition}\label{prop:infinite}
For integers $1\le p<q$ let $I_{p,q}=\bigl(\frac{2p-1}{2q},\frac{2p+1}{2q}\bigr)$
and $f_{p,q}=q\,\mathbf 1_{I_{p,q}}$ on $[0,1]$ with Lebesgue measure $\Leb$.
There is no measurable $E\subseteq[0,1]$ with $\int_Ef_{p,q}\,d\Leb=\frac12$ for
all $p<q$. In particular, the countable family $\{f_{p,q}\}$ has no common
segment.
\end{proposition}

\begin{proof}
Suppose $\Leb(E\cap I_{p,q})=\frac1{2q}$ for all $p<q$. Then $\Leb(E)\ge\Leb(E\cap
I_{1,2})>0$, so by the Lebesgue density theorem~\cite{Fo} there is $x\in
E\cap(0,1)$ with $\Leb(J\setminus E)/\Leb(J)\to0$ for intervals $J$ centred at
$x$ with $\Leb(J)\to0$. For $q$ large, $x$ lies in the closure of some
$I_{p,q}$ with $1\le p<q$ (these closures cover $[\frac1{2q},1-\frac1{2q}]$),
so $I_{p,q}\subseteq J_q=(x-\frac1q,x+\frac1q)$ and
$\Leb(J_q\setminus E)\ge\Leb(I_{p,q}\setminus E)=\frac1{2q}=\frac14\Leb(J_q)$, a
contradiction. If $u$ were a common segment, $E=u(1/2)$ would be such a set.
\end{proof}

\subsection{$\sigma$-finite spaces}

Let $(X,\Sigma,\mu)$ be $\sigma$-finite and non-atomic with $\mu(X)=\infty$. Let
$f$ be measurable with $f-1\in L^1$, and for $E\in\Sigma$ with $\mu(E)<\infty$
put $\mu_f(E)=\mu(E)+\int_E(f-1)\,d\mu$. A \emph{segment} for a family $S$ of
such functions is an increasing map $u\colon[0,\infty)\to\Sigma$ with
$u(0)=\emptyset$, $\bigcup_tu(t)=X$ and $\mu(u(t))=\mu_f(u(t))=t$ for all $t\ge0$
and $f\in S$.

\begin{theorem}\label{thm:sigma}
Let $(X,\Sigma,\mu)$ be $\sigma$-finite and non-atomic with $\mu(X)=\infty$, and
let $S$ be a finite family of measurable functions with $f-1\in L^1$ and
$\int_X(f-1)\,d\mu=0$ for all $f\in S$. Then $S$ admits a segment.
\end{theorem}

\begin{proof}
\emph{Gluing.} Suppose $Y_1\subseteq Y_2\subseteq\cdots$ are measurable,
$\bigcup_kY_k=X$, $\mu(Y_k)<\infty$, and $\mu_f(Y_k)=\mu(Y_k)$ for all $k$ and
$f\in S$. Passing to a subsequence we may assume $0<\mu(Y_1)<\mu(Y_2)<\cdots$;
then $\mu(Y_k)\to\infty$. Put $Y_0=\emptyset$ and $R_k=Y_k\setminus Y_{k-1}$.
Since $\mu_f(R_k)=\mu(R_k)>0$, Theorem~\ref{thm:Fr} applied to $R_k$ yields
$w_k$ with $\mu(w_k(s))=\mu_f(w_k(s))=s\mu(R_k)$ for $f\in S$. Then
$u(t)=Y_{k-1}\cup w_k\bigl((t-\mu(Y_{k-1}))/\mu(R_k)\bigr)$ for
$t\in[\mu(Y_{k-1}),\mu(Y_k)]$ is a segment for $S$.

\emph{Induction on $|S|$.} For $S=\emptyset$ glue along any exhaustion by sets
of finite measure. Suppose segments exist for families of size $n$ and let
$S=\{f_1,\dots,f_{n+1}\}$. Let $v$ be a segment for $\{f_1,\dots,f_n\}$ and
$F(t)=\int_{v(t)}(f_{n+1}-1)\,d\mu$. Since $\mu(v(t)\symd v(s))=|t-s|$, $F$ is
continuous by Lemma~\ref{lem:abscont}; $F(0)=0$, and $F(t)\to0$ as
$t\to\infty$ by dominated convergence. By gluing, it suffices to find sets
$Y_k$ as above for $S$.

\emph{Case 1: $F(t_k)=0$ for a sequence $t_k\to\infty$.} Take $Y_k=v(t_k)$.

\emph{Case 2: the zero set of $F$ is bounded.} Let $T$ be its maximum; we may
assume $F>0$ on $(T,\infty)$, the other case being symmetric. Let $M=\max F>0$.
For $y\in(0,M)$ let $s(y)=\min\{s\ge T:F(s)=y\}$ and $t(y)=\max\{t:F(t)=y\}$.
As $y\downarrow0$, $s(y)$ decreases and $t(y)$ increases, and continuity of $F$
together with $F>0$ on $(T,\infty)$ gives $s(y)\to T$ and $t(y)\to\infty$.
Choose $y_k\downarrow0$, let $s_k=s(y_k)$, $t_k=t(y_k)$, $W=\bigcap_kv(s_k)$ and
\[
Y_k=W\cup\bigl(v(t_k)\setminus v(s_k)\bigr).
\]
The sets $Y_k$ increase and cover $X$: a point outside $W$ lies outside some
$v(s_k)$, hence outside all $v(s_j)$, $j\ge k$, and inside $v(t_j)$ for $j$
large. Moreover $\mu(W)=\mu_{f_i}(W)=T$ for $i\le n$ and $\int_W(f_{n+1}-1)=\lim
F(s_k)=0$, so
$\mu_{f_i}(Y_k)=\mu(Y_k)=T+t_k-s_k$ for $i\le n$, and
$\mu_{f_{n+1}}(Y_k)-\mu(Y_k)=0+F(t_k)-F(s_k)=0$.
\end{proof}

\begin{remark}\label{rem:sigma-counterexample}
Finiteness of $S$ is again necessary. With $f_{p,q}$ as in
Proposition~\ref{prop:infinite}, let $g_{p,q}=f_{p,q}$ on $[0,1]$ and
$g_{p,q}=1$ on $\R\setminus[0,1]$; then $g_{p,q}-1\in L^1(\R)$ has integral $0$.
(The periodic extension of $f_{p,q}$ used in~\cite{Ja} does not satisfy
$g-1\in L^1$.) Suppose $u$ is a segment for all $\mu_{g_{p,q}}$ in the weaker
sense that $\mu_{g_{p,q}}(u(t))=t$ (without requiring a segment for $\mu$). Write
$E_t=u(t)\cap[0,1]$ and $O_t=u(t)\setminus[0,1]$. Then $\int_{E_t}f_{p,q}=t-\Leb(O_t)$
is independent of $(p,q)$; summing over $p$ for fixed $q$ and letting
$q\to\infty$ shows that this common value is $\Leb(E_t)$. Since $\Leb(E_t)$ and
$\Leb(O_t)$ are non-decreasing with continuous sum $t$, both are continuous, and
$\Leb(E_t)$ increases from $0$ to $1$. Choosing $t$ with $\Leb(E_t)=\frac12$
contradicts Proposition~\ref{prop:infinite}.
\end{remark}

%======================================================================
\section{O-segments for a single measure}\label{sec:single}
%======================================================================

We begin with two examples showing what weak outer regularity does and does not
control.

\begin{example}\label{ex:no-O-segment}
Let $X=\bigcup_{n\ge1}I_n$ with $I_n=(2n,2n+1)$, with $\Sigma$ the Borel sets
of $X\subseteq\R$, and $\mu(E)=\sum_n2^{-n}\Leb(E\cap I_n)$. Give $X$ the
topology whose open sets are the unions of subfamilies of $\{I_n\}$. Open sets
are measurable, $\mu$ is non-atomic, and the values of $\mu$ on open sets,
$\{\sum_{n\in N}2^{-n}:N\subseteq\N\}$, fill $[0,1]$. Nevertheless $\mu$ has no
O-segment. Indeed, by Proposition~\ref{prop:usc} an O-segment would give an
upper semicontinuous $\varphi$ with $\varphi_*\mu=\Leb$. Each set $\{\varphi<t\}$
is a union of some $I_n$'s, so $\varphi$ is constant on each $I_n$ (if
$\varphi(x)<t<\varphi(y)$ with $x,y\in I_n$, then $x\in\{\varphi<t\}$ forces
$y\in\{\varphi<t\}$). Hence $\varphi_*\mu$ is purely atomic, a contradiction.
Here $\mu$ is not weakly outer regular: $(2,\tfrac52)$ is not approximable by
open sets.
\end{example}

\begin{example}\label{ex:not-necessary}
Weak outer regularity is sufficient but not necessary. On $X=[0,1]$ with
Lebesgue measure on the Borel sets, take the topology whose open sets are
$\emptyset$, $X$ and the intervals $[0,s)$, $0<s\le1$. Then $u(t)=[0,t)$
($t<1$), $u(1)=X$ is an O-segment, but $\mu$ is not weakly outer regular at
$(\tfrac12,1)$.
\end{example}

\medskip
\noindent\textbf{Standing assumption for Sections~\ref{sec:single} and
\ref{sec:pairs}.} $(X,\Sigma,\mu)$ is a topological measure space and $\mu$ is
finite, non-atomic and weakly outer regular, with $\mu(X)>0$.
\medskip

\begin{lemma}\label{lem:fill}
Let $O\subseteq Y$ be open and $a\in[\mu(O),\mu(Y)]$. There is an open set $O'$
with $O\subseteq O'\subseteq Y$ and $\mu(O')=a$.
\end{lemma}

\begin{proof}
Put $O_0=O$. Given an open $O_n$ with $O\subseteq O_n\subseteq Y$ and
$d_n=a-\mu(O_n)\ge0$, set $O_{n+1}=O_n$ if $d_n=0$. Otherwise, since
$\mu(Y\setminus O_n)\ge d_n$ and $\mu$ is non-atomic, there is a measurable
$A_{n+1}\subseteq Y\setminus O_n$ with $\mu(A_{n+1})=d_n/2$, and by
Lemmas~\ref{lem:approx-equiv} and~\ref{lem:relative-wor} an open
$U_{n+1}\subseteq Y$ with $\mu(A_{n+1}\symd U_{n+1})<d_n/4$. Let
$O_{n+1}=O_n\cup U_{n+1}$. Then
\begin{align*}
\mu(O_{n+1})&\le\mu(O_n)+\mu(U_{n+1})<\mu(O_n)+\tfrac34d_n,\\
\mu(O_{n+1})&\ge\mu(O_n)+\mu(A_{n+1}\cap U_{n+1})>\mu(O_n)+\tfrac14d_n,
\end{align*}
using that $A_{n+1}$ is disjoint from $O_n$. Hence $0\le d_{n+1}\le\frac34d_n$,
so $d_n\to0$, and $O'=\bigcup_nO_n$ is open, $O\subseteq O'\subseteq Y$, and
$\mu(O')=\lim\mu(O_n)=a$.
\end{proof}

\begin{lemma}\label{lem:interpolate}
Let $P\subseteq Q$ be open with $\mu(P)<\mu(Q)$. There is an increasing family
$w(s)$, $s\in[\mu(P),\mu(Q)]$, of open sets with $w(\mu(P))=P$, $w(\mu(Q))=Q$,
$P\subseteq w(s)\subseteq Q$ and $\mu(w(s))=s$.
\end{lemma}

\begin{proof}
Let $\mathcal D$ be the set of dyadic rationals in $[0,1]$ and
$s_r=\mu(P)+r(\mu(Q)-\mu(P))$. We construct open sets $O_r$, $r\in\mathcal D$,
with $O_0=P$, $O_1=Q$, $O_r\subseteq O_{r'}$ for $r<r'$, and $\mu(O_r)=s_r$, by
induction on the dyadic level: if $O_r$ is defined for $r\in2^{-n}\mathbb
Z\cap[0,1]$ and $r=(2k+1)2^{-n-1}$, apply Lemma~\ref{lem:fill} with
$O=O_{k2^{-n}}$, $Y=O_{(k+1)2^{-n}}$ and $a=s_r$. Now let
$w(s)=\bigcup\{O_r:r\in\mathcal D,\ s_r\le s\}$. These sets are open,
increasing in $s$, contained in $Q$, $w(\mu(P))=P$ and $w(\mu(Q))=Q$. Since the
$O_r$ form a chain, $w(s)$ is the union of $O_{r_j}$ along an increasing
sequence $r_j\in\mathcal D$ with $s_{r_j}\uparrow s$, so $\mu(w(s))=s$.
\end{proof}

\begin{proposition}\label{prop:through-open}
For every open set $A$ with $0<\mu(A)<\mu(X)$ there is an O-segment $u$ for
$\mu$ with $u(\mu(A)/\mu(X))=A$.
\end{proposition}

\begin{proof}
Apply Lemma~\ref{lem:interpolate} to the pairs $(\emptyset,A)$ and $(A,X)$ and
concatenate the two families, reparametrized by $t=s/\mu(X)$.
\end{proof}

\begin{proposition}\label{prop:relative-single}
Let $Y$ be open and $A\in\Sigma$ with $A\subseteq Y$ and $\gamma=\mu(Y\setminus
A)>0$. There is an increasing family $v(t)$, $t\in(0,\gamma]$, of open subsets of
$Y$ with $v(\gamma)=Y$, $\mu(A\setminus v(t))=0$ and $\mu(v(t)\setminus A)=t$. If
$A$ is open, one may take $A\subseteq v(t)$ for all $t$.
\end{proposition}

\begin{proof}
By Lemma~\ref{lem:relative-wor} choose open $W_i\subseteq Y$ with
$\mu(A\setminus W_i)=0$ and $\mu(W_i\setminus A)<\gamma2^{-i}$, and put $U_0=Y$,
$U_i=Y\cap W_1\cap\dots\cap W_i$. Then $U_i$ is open and decreasing,
$\mu(A\setminus U_i)=0$, and $s_i=\mu(U_i\setminus A)$ decreases from
$s_0=\gamma$ to $0$. Discarding repetitions, we may assume that either
$s_0>s_1>\cdots$ with $s_i\to0$, or $s_0>\dots>s_N=0$. Since $A\subseteq U_i$
up to a null set, $\mu(V\setminus A)=\mu(V)-\mu(A)$ for every $V\supseteq
U_i$. Apply Lemma~\ref{lem:interpolate} to $U_{i+1}\subseteq U_i$ to get a
family between them, parametrized by $t=\mu(\cdot)-\mu(A)\in[s_{i+1},s_i]$.
These families agree at the common endpoints, and together they define $v$ on
$(0,\gamma]$ with the stated properties. If $A$ is open, replace $v(t)$ by
$v(t)\cup A$; this changes no measure since $\mu(A\setminus v(t))=0$.
\end{proof}

The case $Y=X$ shows that every measurable set $A$ with $\mu(A)<\mu(X)$ is,
up to a null set, the ``initial set'' of a partial O-segment defined on
$(\mu(A)/\mu(X),1]$.

\subsection*{Continuous segments on separable metrizable spaces}

For a single measure on a separable metrizable space, open segments can be
obtained much more directly, and in a stronger form.

\begin{proposition}\label{prop:continuous}
Let $X$ be a separable metrizable space and $\mu$ a finite non-atomic measure on
$\mathcal B(X)$ or on its completion. Then there is a continuous
$\varphi\colon X\to[0,1]$ with $\varphi_*\mu=\mu(X)\Leb$. Consequently
$u(t)=\{\varphi<t\}$ $(t<1)$, $u(1)=X$, is an O-segment for $\mu$, and the
closed sets $K(t)=\{\varphi\le t\}$ satisfy $u(t)\subseteq K(t)$ and
$\mu(K(t)\setminus u(t))=0$.
\end{proposition}

\begin{proof}
Fix a compatible metric $d\le1$ and a dense sequence $(x_k)$, and put
$f_k(x)=d(x,x_k)$; these continuous functions separate points. Let $P$ be the
product of Lebesgue measures on $\Omega=[0,1]^{\N}$, and for
$c=(c_k)\in\Omega$ let $\varphi_c(x)=\sum_k2^{-k}c_kf_k(x)$. The map
$(c,x)\mapsto\varphi_c(x)$ is continuous on $\Omega\times X$, so
$E=\{(c,x,y):\varphi_c(x)=\varphi_c(y)\}$ is closed in the separable metrizable
space $\Omega\times X\times X$, whose Borel $\sigma$-algebra is the product
$\sigma$-algebra. For $x\ne y$ choose $k$ with $f_k(x)\neq f_k(y)$; for fixed
values of the other coordinates, $\varphi_c(x)-\varphi_c(y)$ is a non-constant
affine function of $c_k$, so $P(\{c:\varphi_c(x)=\varphi_c(y)\})=0$. Since $\mu$
is non-atomic and singletons are Borel, $\mu(\{x\})=0$ for all $x$, so the
diagonal of $X\times X$ is $\mu\otimes\mu$-null. By Tonelli,
\[
\int_\Omega(\mu\otimes\mu)\bigl(\{(x,y):\varphi_c(x)=\varphi_c(y)\}\bigr)\,dP(c)
=\int_{X\times X}P\bigl(\{c:\varphi_c(x)=\varphi_c(y)\}\bigr)\,d(\mu\otimes\mu)=0 .
\]
Fix $c$ for which the integrand vanishes. If $\nu=(\varphi_c)_*\mu$ had an atom
$\{r\}$ of mass $m>0$, the integrand would be at least $m^2$; hence $\nu$ has a
continuous distribution function $\Phi(r)=\nu((-\infty,r])/\mu(X)$. Then
$\varphi=\Phi\circ\varphi_c$ is continuous with values in $[0,1]$, and for
$t\in(0,1)$, with $r_t=\min\{r:\Phi(r)=t\}$, we get
$\mu(\{\varphi<t\})=\nu((-\infty,r_t))=t\mu(X)$. Thus $\varphi_*\mu=\mu(X)\Leb$,
and the remaining statements follow from Proposition~\ref{prop:usc}(b) and
$\mu(\{\varphi=t\})=0$.
\end{proof}

\begin{remark}\label{rem:support}
(1) The same conclusion holds on an arbitrary metrizable space whenever $\mu$
is concentrated on a separable closed set $S$ (e.g.\ if $\mu$ is a Radon
measure): apply Proposition~\ref{prop:continuous} on $S$ and extend $\varphi$
to $X$ by the Tietze extension theorem.

(2) If $X=[0,1]^n$ and $\mu$ is non-atomic, of full support and
$\mu(\partial X)=0$, the homeomorphic measures theorem of Oxtoby and
Ulam~\cite{OU,AP} gives a homeomorphism $h$ of $X$ with
$h_*\mu=\mu(X)\Leb^n$, and the first coordinate of $h$ is such a $\varphi$.

(3) Proposition~\ref{prop:continuous} does not extend in an obvious way to two
densities: the randomization argument controls atoms of a single pushforward,
not the joint distribution of $\mu$ and $\mu_f$. This is the situation treated
in Section~\ref{sec:pairs}.
\end{remark}

%======================================================================
\section{Common O-segments for pairs of functions}\label{sec:pairs}
%======================================================================

The standing assumption of Section~\ref{sec:single} remains in force. The idea
of the proof is to construct the O-segment along a dense set of parameters,
each new open set being obtained by solving one scalar equation with the
intermediate value theorem. That equation is solved along a one-parameter
family of open sets on which a second measure serves as the parameter; the
following definition isolates what is needed of such families.

\begin{definition}\label{def:path}
Let $g\in L^1$, let $Y$ be open and let $A\in\Sigma$ with $A\subseteq Y$ and
$\gamma=\mu_g(Y\setminus A)>0$. A \emph{$g$-path from $A$ to $Y$} is an increasing
family $v(t)$, $t\in(0,\gamma]$, of open subsets of $Y$ such that $v(\gamma)=Y$,
$\mu(A\setminus v(t))=0$, $\mu_g(v(t)\setminus A)=t$ for all $t$, and the
function $t\mapsto\mu(v(t)\setminus A)$ is continuous on $(0,\gamma]$ and tends to
$0$ as $t\downarrow0$. We say that $g$ has \emph{property (P)} if a $g$-path
from $A$ to $Y$ exists for all such $Y$ and $A$.
\end{definition}

\begin{remark}\label{rem:path-continuity}
If $v$ is a $g$-path from $A$ to $Y$ and $s<t$, then $v(t)\symd
v(s)=v(t)\setminus v(s)$ coincides up to a null set with $(v(t)\setminus
A)\setminus(v(s)\setminus A)$, so $\mu(v(t)\symd v(s))=\mu(v(t)\setminus
A)-\mu(v(s)\setminus A)$; similarly $\mu(v(t)\symd A)=\mu(v(t)\setminus A)$.
Consequently, by Lemma~\ref{lem:abscont}, for every $k\in L^1$ the function
$t\mapsto\mu_k(v(t)\setminus A)$ is continuous on $(0,\gamma]$ and tends to $0$
as $t\downarrow0$.
\end{remark}

\begin{lemma}\label{lem:constants}
Every positive constant $\kappa$ has property (P).
\end{lemma}

\begin{proof}
Given $Y$ and $A$ with $\kappa\mu(Y\setminus A)=\gamma$, take the family $w$ of
Proposition~\ref{prop:relative-single} and put $v(t)=w(t/\kappa)$.
\end{proof}

\begin{lemma}[Key lemma]\label{lem:key}
Let $g\in L^1$ have property (P) and let $f\in L^1$. Let $Y$ be open and
$B\in\Sigma$ with $B\subseteq Y$ and
\[
\mu_f(Y\setminus B)=\mu_g(Y\setminus B)=c>0 .
\]
Then for every $a\in[0,c)$ and $\eps\in(0,c-a]$ there is an open set
$U\subseteq Y$ with $\mu(B\setminus U)=0$ and
\[
\mu_f(U\setminus B)=\mu_g(U\setminus B)\in(a,a+\eps).
\]
\end{lemma}

\begin{proof}
Put $h=f-g$.

\emph{Step 1: an intermediate measurable set.} Since $\mu_g(Y\setminus
B)\ne0$, we have $\mu(Y\setminus B)>0$, and Theorem~\ref{thm:Fr} applied to
$Y\setminus B$ with $S=\{f,g\}$ yields a segment $w$ on $Y\setminus B$ with
$\mu_f(w(s))=\mu_g(w(s))=sc$. Let $A=B\cup w(a/c)$. Then $B\subseteq
A\subseteq Y$, $\mu_f(A\setminus B)=\mu_g(A\setminus B)=a$, $\mu_h(A\setminus
B)=0$ and $\mu_g(Y\setminus A)=c-a>0$.

\emph{Step 2: a path and the function $F$.} Let $v$ be a $g$-path from $A$ to
$Y$, and for $s\in(a,c]$ put $V_s=v(s-a)$; put $V_a=A$. Then $V_s$ is open for
$s>a$, $V_s\subseteq V_{s'}$ for $a<s\le s'$, $A\subseteq V_s$ up to a null set,
$V_c=Y$, and
\[
\mu_g(V_s\setminus B)=\mu_g(A\setminus B)+\mu_g(V_s\setminus A)=s .
\]
Define $F(s)=\mu_h(V_s\setminus B)=\mu_h(V_s\setminus A)$ for $s\in[a,c]$. By
Remark~\ref{rem:path-continuity}, $F$ is continuous on $[a,c]$; moreover
$F(a)=0$ and $F(c)=\mu_h(Y\setminus B)=0$. More generally, for $k\in L^1$,
$s\mapsto\mu_k(V_s\setminus B)$ is continuous on $[a,c]$.

If $F(s)=0$ for some $s\in(a,a+\eps)$, then $U=V_s$ has the required
properties: $\mu(B\setminus U)\le\mu(A\setminus V_s)=0$ and
$\mu_f(U\setminus B)=\mu_g(U\setminus B)=s$. Assume from now on that $F$ has no
zero in $(a,a+\eps)$. Let $t_0=\min\{s\in[a+\eps,c]:F(s)=0\}$, which exists
since $F(c)=0$. Then $F$ has no zero in $(a,t_0)$, so it has constant sign
there.

\emph{Step 3: choice of parameters.} Fix $\tau\in(a,a+\eps/4)$. Since
$F(s)\to0$ as $s\uparrow t_0$, we can choose $t_2\in(\max\{\tau,t_0-\eps/4\},t_0)$
with $0<|F(t_2)|<|F(\tau)|/2$. Then $2F(t_2)$ lies strictly between $F(a)=0$ and
$F(\tau)$, so there is $t_3\in(a,\tau)$ with $F(t_3)=2F(t_2)$. Thus
\[
a<t_3<\tau<t_2<t_0,\qquad t_3-a<\eps/4,\qquad t_0-t_2<\eps/4 .
\]

\emph{Step 4: an open correction set.} Let $C=V_{t_0}\setminus V_{t_2}$. Since
$B\subseteq A\subseteq V_{t_2}$ up to null sets, $C\cap B$ is null, and
\[
\kappa:=\mu_h(C)=F(t_0)-F(t_2)=-F(t_2)\neq0,\qquad \mu_g(C)=t_0-t_2>0 .
\]
Let $\theta=\min\{|\kappa|/4,\,(t_0-t_2)/2,\,\eps/4\}$, choose $\delta>0$ by
Lemma~\ref{lem:abscont} so that $\int_E(|f|+|g|)\,d\mu<\theta$ whenever
$\mu(E)<\delta$, and choose by Lemma~\ref{lem:relative-wor} an open $W\subseteq
Y$ with $\mu(C\setminus W)=0$ and $\mu(W\setminus C)<\delta$. Put $R=W\setminus
C$.

\emph{Step 5: the intermediate value argument.} For $s\in[a,t_3]$ let
$G(s)=\mu_h((V_s\cup W)\setminus B)$. As in Step~2, $G$ is continuous on
$[a,t_3]$, because $\mu\bigl((V_s\cup W)\symd(V_{s'}\cup W)\bigr)\le\mu(V_s\symd
V_{s'})$. For $s\le t_3$ we have $V_s\subseteq V_{t_2}$ up to a null set, so
$V_s\cap C$ is null; also $C\subseteq W$ up to a null set and $C\cap B$ is
null. Hence, up to null sets, $(V_s\cup W)\setminus B$ is the disjoint union of
$V_s\setminus B$, $C$ and $R_s=W\setminus(C\cup V_s\cup B)\subseteq R$, and
\[
G(s)=F(s)+\kappa+e(s),\qquad |e(s)|\le\int_R|h|\,d\mu<\theta\le|\kappa|/4 .
\]
Therefore $G(a)=\kappa+e(a)$ has the sign of $\kappa$, while
$G(t_3)=2F(t_2)+\kappa+e(t_3)=-\kappa+e(t_3)$ has the sign of $-\kappa$. By the
intermediate value theorem there is $s^*\in(a,t_3)$ with $G(s^*)=0$.

\emph{Step 6: conclusion.} Let $U=V_{s^*}\cup W$. Then $U$ is open,
$U\subseteq Y$, $\mu(B\setminus U)\le\mu(A\setminus V_{s^*})=0$, and
$\mu_f(U\setminus B)-\mu_g(U\setminus B)=G(s^*)=0$. By the same decomposition,
\[
\mu_g(U\setminus B)=s^*+(t_0-t_2)+e',\qquad |e'|\le\int_R|g|\,d\mu<\theta .
\]
Since $\theta\le(t_0-t_2)/2$, we get $\mu_g(U\setminus B)>s^*>a$, and since
$\theta\le\eps/4$,
$\mu_g(U\setminus B)<t_3+(t_0-t_2)+\eps/4<a+\tfrac34\eps$.
\end{proof}

\begin{proposition}\label{prop:relative-P}
Let $g\in L^1$ have property (P) and let $f\in L^1$. Let $Y$ be open and
$D\in\Sigma$ with $D\subseteq Y$ and $\mu_f(Y\setminus D)=\mu_g(Y\setminus
D)=\gamma>0$. Then there is an increasing family $u(t)$, $t\in(0,\gamma]$, of
open subsets of $Y$ such that $u(\gamma)=Y$ and, for all $t$,
\[
\mu(D\setminus u(t))=0,\qquad \mu_f(u(t)\setminus D)=\mu_g(u(t)\setminus D)=t .
\]
If $D$ is open, one may take $D\subseteq u(t)$ for all $t$.
\end{proposition}

\begin{proof}
We construct finite sets $\{\gamma\}=Q_0\subseteq Q_1\subseteq\cdots\subseteq
(0,\gamma]$ and open sets $O_q\subseteq Y$, $q\in Q=\bigcup_nQ_n$, such that
\begin{equation}\label{eq:chain}
O_\gamma=Y,\quad O_q\subseteq O_{q'}\ (q<q'),\quad \mu(D\setminus O_q)=0,\quad
\mu_f(O_q\setminus D)=\mu_g(O_q\setminus D)=q,
\end{equation}
and such that every gap of $\{0\}\cup Q_n$ has length at most
$(2/3)^n\gamma$; if $D$ is open we also require $D\subseteq O_q$.

Suppose $Q_n$ and $O_q$, $q\in Q_n$, are constructed. Let $q'<q''$ be
consecutive elements of $\{0\}\cup Q_n$. Put $Y'=O_{q''}$, and put $B=O_{q'}$ if
$q'>0$ and $B=D\cap Y'$ if $q'=0$. In both cases $B\subseteq Y'$, $\mu(D\setminus B)=0$
by~\eqref{eq:chain}, and $\mu_f(B\setminus D)=\mu_g(B\setminus D)=q'$. Since
$(Y'\cap D)\setminus B$ is null, it follows that
\[
\mu_f(Y'\setminus B)=\mu_g(Y'\setminus B)=q''-q'=:c>0 .
\]
Apply Lemma~\ref{lem:key} with $a=\eps=c/3$ to obtain an open $U\subseteq Y'$
with $\mu(B\setminus U)=0$ and $\mu_f(U\setminus B)=\mu_g(U\setminus
B)=:\beta\in(c/3,2c/3)$. Let $O_q=U\cup B$ if $q'>0$, and $O_q=U$ (or $U\cup D$
if $D$ is open) if $q'=0$, where $q=q'+\beta$. Then $O_q$ is open,
$O_{q'}\subseteq O_q\subseteq O_{q''}$ (for $q'=0$ only the second inclusion is
needed), $\mu(D\setminus O_q)\le\mu(D\setminus B)+\mu(B\setminus U)=0$, and
\[
\mu_f(O_q\setminus D)=\mu_f(B\setminus D)+\mu_f(U\setminus B)=q'+\beta=q,
\]
and likewise for $g$; here we used that $(U\cap D)\setminus B$ is null and that
modifying a set by a null set does not change $\mu_f$ or $\mu_g$. Doing this for
every gap gives $Q_{n+1}$; the two new gaps have length less than $2c/3$.

The set $Q$ is dense in $[0,\gamma]$. For $t\in(0,\gamma]$ let
$u(t)=\bigcup\{O_q:q\in Q,\ q\le t\}$. Then $u$ is increasing, $u(t)$ is an
open subset of $Y$, $u(\gamma)\supseteq O_\gamma=Y$, and $\mu(D\setminus u(t))=0$
since $Q\cap(0,t]\neq\emptyset$. As the $O_q$ form a chain, $u(t)$ is the union
of $O_{q_j}$ along an increasing sequence $q_j\in Q$ with $q_j\uparrow t$, so
$\mu_f(u(t)\setminus D)=\lim q_j=t$, and likewise for $g$. If $D$ is open, all
$O_q$ contain $D$, hence so does $u(t)$.
\end{proof}

\begin{corollary}\label{cor:allP}
Every $g\in L^1$ has property (P).
\end{corollary}

\begin{proof}
Let $Y$ be open and $A\subseteq Y$ with $\gamma=\mu_g(Y\setminus A)>0$. Then
$\mu(Y\setminus A)>0$; put $\kappa=\gamma/\mu(Y\setminus A)>0$. The constant
$\kappa$ has property (P) by Lemma~\ref{lem:constants}, and
$\mu_g(Y\setminus A)=\mu_\kappa(Y\setminus A)=\gamma$. Proposition~\ref{prop:relative-P},
applied to the pair $(g,\kappa)$ with $D=A$, gives an increasing family $v$ of
open subsets of $Y$ with $v(\gamma)=Y$, $\mu(A\setminus v(t))=0$ and
$\mu_g(v(t)\setminus A)=\kappa\mu(v(t)\setminus A)=t$. Thus $\mu(v(t)\setminus
A)=t/\kappa$ is continuous and tends to $0$, and $v$ is a $g$-path from $A$ to
$Y$.
\end{proof}

Combining Proposition~\ref{prop:relative-P} with Corollary~\ref{cor:allP} we
obtain the main technical result.

\begin{theorem}\label{thm:relative}
Let $(X,\Sigma,\mu)$ be a topological measure space with $\mu$ finite,
non-atomic and weakly outer regular, and let $f,g\in L^1$. Let $Y$ be open and
$D\in\Sigma$ with $D\subseteq Y$ and $\mu_f(Y\setminus D)=\mu_g(Y\setminus
D)=\gamma>0$. Then there is an increasing family $u(t)$, $t\in(0,\gamma]$, of
open subsets of $Y$ with $u(\gamma)=Y$ such that for all $t$
\[
\mu(D\setminus u(t))=0\quad\text{and}\quad\mu_f(u(t)\setminus D)=\mu_g(u(t)\setminus D)=t .
\]
If $D$ is open, one may take $D\subseteq u(t)$ for all $t$.
\end{theorem}

\begin{theorem}\label{thm:main}
Let $(X,\Sigma,\mu)$ be a topological measure space with $\mu$ finite,
non-atomic and weakly outer regular.
\begin{enumerate}
\item[(a)] For every $f\in L^1$ there is an O-segment for $\mu$ which is also a
segment for $f$.
\item[(b)] For all $f,g\in L^1$ with $(\mu_f(X),\mu_g(X))\neq(0,0)$ there is a
common O-segment for $f$ and $g$.
\item[(c)] In (b), let in addition $D$ be an open set and $s\in(0,1)$ with
$\mu_f(D)=s\mu_f(X)$ and $\mu_g(D)=s\mu_g(X)$. Then the common O-segment can be
chosen with $u(s)=D$.
\end{enumerate}
\end{theorem}

\begin{proof}
(b) Let $\alpha=\mu_f(X)$ and $\beta=\mu_g(X)$. Choose $f',g'$ in the span of
$f,g$ with $\mu_{f'}(X)=\mu_{g'}(X)=1$ such that $f,g$ lie in the span of
$f',g'$: if $\alpha\beta\ne0$, take $f'=f/\alpha$, $g'=g/\beta$; if $\alpha=0\ne\beta$,
take $g'=g/\beta$, $f'=f+g'$ (so that $f=f'-g'$); the case $\beta=0\ne\alpha$ is
symmetric. Theorem~\ref{thm:relative} with $Y=X$, $D=\emptyset$, $\gamma=1$ and
the pair $(f',g')$ gives $u$ on $(0,1]$ with $u(1)=X$ and
$\mu_{f'}(u(t))=\mu_{g'}(u(t))=t$; setting $u(0)=\emptyset$ we obtain a common
O-segment for $f'$ and $g'$, hence for $f$ and $g$ by Remark~\ref{rem:linear}(1).

(a) This is (b) with $g\equiv1$, since $\mu_1(X)=\mu(X)>0$.

(c) With $f',g'$ as in (b) we have $\mu_{f'}(D)=\mu_{g'}(D)=s$ by linearity.
Theorem~\ref{thm:relative} with $Y=D$ and the initial set $\emptyset$ gives a
family $u_1$ on $(0,s]$ with $u_1(s)=D$ and $\mu_{f'}(u_1(t))=\mu_{g'}(u_1(t))=t$.
Theorem~\ref{thm:relative} with $Y=X$ and the open initial set $D$ gives a
family $u_2$ on $(0,1-s]$ with $D\subseteq u_2(t)$, $u_2(1-s)=X$ and
$\mu_{f'}(u_2(t)\setminus D)=\mu_{g'}(u_2(t)\setminus D)=t$. Put $u(0)=\emptyset$,
$u(t)=u_1(t)$ for $0<t\le s$ and $u(t)=u_2(t-s)$ for $s<t\le1$.
\end{proof}

\begin{remark}\label{rem:special-cases}
(1) Theorem~\ref{thm:relative} with $f=g=1$ recovers
Proposition~\ref{prop:relative-single}, and Theorem~\ref{thm:main}(c) with
$f=g=1$ recovers Proposition~\ref{prop:through-open}.

(2) If $g\ge0$ and $\{g>0\}$ has full measure, Theorem~\ref{thm:main}(b) can
also be deduced from (a) applied to the measure $\nu=g\,d\mu$, which is
non-atomic and, being equivalent to $\mu$, weakly outer regular, and to the
function $f/g\in L^1(\nu)$. The point of the proof above is that no positivity
of $g$ is needed: a $g$-path is built from the case of constant densities, and
Lemma~\ref{lem:key} only uses its parametrization by $\mu_g$ together with the
continuity in $\mu$ supplied by Definition~\ref{def:path}.
\end{remark}

%======================================================================

\end{document}